\documentclass[aop,MSNbibl,nameyear,dvips]{arximspdf}
\usepackage{mathbh}

\doi{10.1214/12-AOP819} 
\volume{41}
\issue{5}
\pubyear{2013}
\firstpage{3157}
\lastpage{3180}

\makeatletter

\newcommand{\iiint}{\int\!\!\int\!\!\int}

\newproclaim{remark}{Remark}

\newproclaim{assumption}{Assumption}

\newcommand{\rrVert}{\Vert}
\newcommand{\rrvert}{\vert}
\newcommand{\llVert}{\Vert}
\newcommand{\llvert}{\vert}

\newtheorem{theorem}{Theorem}[section]
\newtheorem{lemma}[theorem]{Lemma}
\newtheorem{proposition}[theorem]{Proposition}

\newproclaim{definition}[theorem]{Definition}

\newcommand{\interleave}{\Vert\hspace*{-1.4pt}\vert}

\newcommand{\interleavee}{\big\Vert\hspace*{-2pt}\big\vert}

\newcommand{\underset}[2]{\mathop{#2}_{#1}}

\newcommand{\Var}{\operatorname{Var}}

\newcommand{\bE}{\mathbb{E}}
\newcommand{\bN}{\mathbb{N}}
\newcommand{\bP}{\mathbb{P}}
\newcommand{\bR}{\mathbb{R}}
\newcommand{\bZ}{\mathbb{Z}}

\newcommand{\ind}{\mathbh{1}}

\newcommand{\cC}{\mathcal{C}}
\newcommand{\cF}{\mathcal{F}}

\newcommand{\fF}{\mathfrak{F}}

\newcommand{\bs}{\setminus}

\newcommand{\konv}[1]{\underset{#1}{\longrightarrow}}

\newcommand{\erwsymbol}{\mathbb{E}}

\makeatother

\begin{document}
\begin{frontmatter}

\title{Random walks in dynamic random environments: A~transference principle\thanksref{T1}}
\runtitle{RWDRE: A transference principle}

\thankstext{T1}{Supported by the NWO, under the project ``vrije competitie'' number 600.065.100.\break07N14.}

\begin{aug}
\author[A]{\fnms{Frank} \snm{Redig}\ead[label=e1]{F.H.J.Redig@tudelft.nl}}
\and
\author[B]{\fnms{Florian} \snm{V\"ollering}\corref{}\ead[label=e2]{florian.voellering@mathematik.uni-goettingen.de}}
\runauthor{F. Redig and F. V\"ollering}
\affiliation{Delft University of Technology and
Georg-August-Universit\"at G\"ottingen}
\address[A]{Institute of Mathematics\\
Delft University of Technology\\
Mekelweg 4\\
2628 CD Delft\\
The Netherlands\\
\printead{e1}}
\address[B]{Institut f\"ur Mathematische Stochastik\\
Georg-August-Universit\"at G\"ottingen\\
Goldschmidtstra{\ss}e 7\\
37077 G\"ottingen\\
Germany\\
\printead{e2}} 
\end{aug}

\received{\smonth{7} \syear{2011}}
\revised{\smonth{8} \syear{2012}}

%
\begin{abstract}
We study a general class of random walks driven by a uniquely ergodic
Markovian environment. Under a coupling condition on the environment we
obtain strong ergodicity properties for the environment as seen from
the position of the walker, that is, the environment process. We can
transfer the rate of mixing in time of the environment to the rate of
mixing of the environment process with a loss of at most polynomial
order. Therefore the method is applicable to environments with
sufficiently fast polynomial mixing. We obtain unique ergodicity of the
environment process. Moreover, the unique invariant measure of the
environment process depends continuously on the jump rates of the
walker.

As a consequence we obtain the law of large numbers and a central limit
theorem with nondegenerate variance for the position of the walk.
\end{abstract}

%
\begin{keyword}[class=AMS]
\kwd[Primary ]{82C41}
\kwd[; secondary ]{60F17}
\end{keyword}
\begin{keyword}
\kwd{Environment process}
\kwd{coupling}
\kwd{random walk}
\kwd{transference principle}
\kwd{central limit theorem}
\end{keyword}

\end{frontmatter}

\section{Introduction}
In recent days random walks in dynamic random environments have been studied
by several authors. Motivation comes among others from nonequilibrium
statistical
mechanics---derivation of Fourier law---[\citet{DOLGOPYATLIVERANI08}]
and large deviation theory
[\citet{RASSOULAGHASAPPALAINEN11}].
In principle, random walks in dynamic random environments contain, as a
particular case, a
random walk in a static random environment. However, mostly, in turning
to dynamic
environments, authors concentrate more on environments with sufficient
mixing properties. In that case the fact that the environment
is dynamic helps to obtain self-averaging properties that ensure
standard limiting behavior of the walk, that is, the law of large numbers
and the central limit theorem.

In the study of the limiting behavior of the walker, the environment
process, that is, the environment as seen from the position
of the walker plays a crucial role. See also \citet
{JOSEPHRASSOULAGHA10}, \citet{RASSOULAGHA03} for the use of the
environment process in related context. In a translation invariant
setting, the environment process is a Markov
process and its ergodic properties fully determine corresponding
ergodic properties of the walk, since the position of the walker
equals an additive function of the environment process plus a
controllable martingale.

The main theme of this paper is the following
natural question: if the environment is uniquely ergodic, with
a sufficient speed of mixing, then the environment process
shares similar properties.
In several works [\citet{BOLDREGHINIIGNATYUK92}, \citet
{BANDYOPADHYAYZEITOUNI06},
\citet{BOLDREGHINIMINLOS07}, \citet{AVENADENHOLLANDERREDIG11}] this transfer
of ``good properties of the environment'' to ``similar properties
of the environment process'' is made via a perturbative argument,
and therefore holds only in a regime where the environment
and the walker are weakly coupled. Some nonperturbative results also
exist, but those require strong mixing properties of the environment in space
and time [\citet{DOLGOPYATKELLERLIVERANI08}, \citet
{DOLGOPYATLIVERANI09}, \citet{BRICMONTKUPIAINEN09}].

In this paper we consider the context of general Markovian uniquely
ergodic environments, which are such that the semigroup contracts at a
minimal speed in a norm of variation type. Examples of such
environments include interacting particle systems in ``the $M<\varepsilon
$ regime'' [\citet{LIGGETT05}] and weakly interacting diffusion
processes on a compact manifold. Our conditions on the environment are
formulated in the language of coupling. More precisely, we impose that
for the environment there exists a coupling such that the distance
between every pair of
initial configurations in this coupling decays fast enough so that
multiplied with $t^d$ it is still integrable in time. As a result we
then obtain that for the environment process there exists a coupling
such that the distance between every pair of initial configurations in
this coupling decays at a speed which is at least
integrable in time. In fact we show more, namely in going from the
environment to the environment process, we essentially loose a factor
$t^d$ in the rate of decay to equilibrium. For example, if for the
environment there is a coupling where the distance decays
exponentially, then the same holds for the environment process (with
possibly another rate).

Once we have controllable coupling properties of the environment
process, we can draw strong conclusions for the position of the walker,
for example, a law of large numbers with an asymptotic speed that
depends continuously on the rates, and a central limit
theorem. We also prove recurrence in $d=1$ under condition of zero speed.

Our paper is organized as follows.
The model and necessary notation are introduced in Section \ref
{sectionmodel}. Section~\ref{sectionEP} is dedicated to lift
properties of the environment to the environment process. The focus is
on Theorem~\ref{thmEP-estimate} and its refinements. Based on these
results consequences for the walker are summarized in Section \ref
{sectionconsequences}. In Section~\ref{sectionexamples} we give
examples for environments to which the results
are applicable and present one example which has polynomial mixing in
space and time. Section~\ref{sectionproofs} is devoted to proofs.

\section{The model}\label{sectionmodel}
\subsection{Environment}
A random walk in dynamic random environment is a process $(X_t)_{t\geq
0}$ on the lattice $\bZ^d$ which is driven by a second process $(\eta
_t)_{t\geq0}$ on $E^{\bZ^d}$, the (dynamic) environment. This is
interpreted as a random walk moving through the environment, with
time-dependent transition rates being determined by the local
environment around the random walk.

To become more precise, the \textit{environment} $(\eta_t)_{t\geq0}$
we assume to be a Feller process on the state space $\Omega:=E^{\bZ
^d}$, where $(E,\rho)$ is a compact metric space (examples in mind are
$E=\{0,1\}$ or $E=[0,1]$). We assume (without loss of generality) that
the distance $\rho$ on $E$ is bounded from above by 1. The generator
of the Markov process $(\eta_t)_{t\geq0}$ is denoted by $L^E$ and its
semigroup by $S_t^E$, both considered on the space of continuous
functions $\cC(\Omega;\bR)$. We assume that the environment is
translation invariant, that is,
\[
\bP^E_\eta(\theta_x\eta_t \in
\cdot) = \bP^E_{\theta_x\eta
}(\eta_t \in\cdot)
\]
with $\theta_x$ denoting the shift operator $\theta_x\eta(y) = \eta
(y-x)$ and $\bP^ E_\eta$ the path space measure of the process $(\eta
_t)_{t\geq0}$ starting from $\eta$.

\subsection{Lipschitz functions}

Denote, for $x\in\bZ^d$,
\begin{eqnarray}
&&(\Omega\times\Omega)_{x}:= \bigl\{(\eta,\xi) \in\Omega^2\dvtx
\eta(x) \neq\xi (x) \mbox{ and } \eta(y)=\xi(y) \ \forall y\in\bZ^d
\bs\{x\} \bigr\},\nonumber\\
&&\eqntext{x\in\bZ^d.}
\end{eqnarray}

\begin{definition}
For any $f\dvtx \Omega\to\bR$, we denote by $\delta_f(x)$ the
Lipschitz-constant of $f$ w.r.t. the variable $\eta(x)$,
\[
\delta_f(x):= \sup_{(\eta,\xi)\in(\Omega\times\Omega
)_{x}}\frac{f(\eta
)-f(\xi)}{\rho(\eta(x),\xi(x))}.
\]
We write
%
\begin{equation}
\interleave f \interleave:=\sum_{x\in\bZ^d}
\delta_f(x).
\end{equation}
\end{definition}
Note that $\interleave f \interleave<\infty$ implies that $f$ is
bounded and
continuous, and the value of $f$ is uniformly weakly dependent on sites
far away. A weaker semi-norm we also use is the oscillation (semi)-norm
\[
\llVert f \rrVert _{\mathrm{osc}}:= \sup_{\eta,\xi\in\Omega} \bigl(f(
\eta)-f(\xi ) \bigr).
\]
From telescoping over single site changes one sees $\llVert  f
\rrVert _{\mathrm{osc}}
\leq\interleave f \interleave$.

\subsection{The random walker and assumption on rates}\label{subsectionrates}

The random walk $X_t$ is a process on $\bZ^d$, whose transition rates
depend on the state of the environment as seen from the walker. More
precisely, the rate to jump from site $x$ to site $x+z$ given that the
environment is in state $\eta$ is $\alpha(\theta_{-x}\eta,z)$.
We make two assumptions on these jump rates $\alpha$. First, we
guarantee that the position of the walker $X_t$ has a first moment by assuming
%
\begin{equation}
\label{eqrates-moment} \llVert \alpha\rrVert _1:= \sum
_{z\in\bZ^d}\llVert z \rrVert \sup_{\eta\in
\Omega}\bigl
\llvert \alpha(\eta,z) \bigr\rrvert <\infty.
\end{equation}
More generally, as sometimes higher moments are necessary, we write
\[
\llVert \alpha\rrVert _p^p:= \sum
_{z\in\bZ^d}\llVert z \rrVert ^p\sup_{\eta\in
\Omega}
\bigl\llvert \alpha(\eta,z) \bigr\rrvert,\qquad p\geq1.
\]
Second, we limit the sensitivity of the rates to small changes in the
environment by assuming that
%
\begin{equation}
\interleave\alpha\interleave:=\sum_{z\in\bZ^d} \interleavee
\alpha(\cdot,z) \interleavee<\infty.
\end{equation}
Finally, sometimes we will have to assume the stronger estimate
%
\begin{equation}
\interleave\alpha\interleave_1:=\sum_{z\in\bZ^d}
\llvert z \rrvert \interleavee\alpha(\cdot,z) \interleavee<\infty.
\end{equation}

\subsection{Environment process}
While the random walker $X_t$ itself is not a Markov process due to the
dependence on the environment, the pair $(\eta_t,X_t)$ is a Markov
process with generator
\[
Lf(\eta,x) = L^Ef(\cdot,x) (\eta) + \sum
_{z\in\bZ^d}\alpha (\theta_{-x}\eta,z) \bigl[f(\eta,x+z)-f(
\eta,x) \bigr],
\]
corresponding semigroup $S_t$ (considered on the space of functions
continuous in $\eta\in\Omega$ and Lipschitz continuous in $x\in\bZ
^d$) and path space measure $\bP_{\eta,x}$.

The environment as seen from the walker is of crucial importance to
understand the asymptotic behavior of the walker itself. This
process,\break
$(\theta_{-X_t}\eta_t)_{t\geq0}$, is called the \textit{environment
process} (not to be confused with the environment $\eta_t$). It is a
Markov process with generator
\[
L^{\mathrm{EP}}f(\eta) = L^E f(\eta) + \sum
_{z\in\bZ^d}\alpha(\eta,z) \bigl[f(\theta_{-z}\eta)-f(\eta)
\bigr],
\]
corresponding semigroup $S_t^{\mathrm{EP}}$ [on $\cC(\Omega)$] and path space
measure $\bP^{\mathrm{EP}}_{\eta}$. Notice that this process is meaningful
only in the translation invariant context.

\subsection{Coupling of the environment}

In the remainder\vspace*{1pt} of the paper we will need a coupling
$\widehat{\bP}{}^E_{\eta,\xi}$, $\eta,\xi\in\Omega$, of the environment.
For $\eta,\xi$ the coupled pair $(\eta^1_t,\eta^2_t)_{t\geq0}$ consists
of two copies of the environment, started\vspace*{1pt} in $\eta$ and $\xi$. By
definition, a~coupling has the marginals\vspace*{2pt}
$\widehat{\bP}{}^E_{\eta,\xi}(\eta^1\in\cdot)=\bP_\eta(\eta_t\in\cdot)$
and $\widehat{\bP}{}^E_{\eta,\xi
}(\eta^2\in\cdot)=\bP_\xi(\eta_t\in\cdot)$. Let $(\cF_t)_{t\geq0}$ be
the canonical filtration in the path space of coupled processes. We say
such a coupling satisfies the marginal Markov property if, for any
$f\dvtx \Omega\to\bR$,
%
\begin{equation}
\label{eqMMP} \widehat{\bE}{}^E_{\eta,\xi} \bigl[ f\bigl(
\eta^i_{t}\bigr) | \cF_s \bigr] =
S_{t-s}^E f\bigl(\eta^i_s\bigr),\qquad
i=1,2; t\geq s \geq0.
\end{equation}
We say it satisfies the strong marginal Markov property if, for any
$f\dvtx \Omega\to\bR$ and any stopping time $\tau$,
%
\begin{equation}
\label{eqSMMP} \widehat{\bE}{}^E_{\eta,\xi} \bigl[
\ind_{t\geq\tau}f\bigl(\eta^i_t\bigr) |
\cF_\tau \bigr] = \ind_{t\geq\tau} S_{t-\tau}^E f
\bigl(\eta ^i_\tau\bigr),\qquad i=1,2.
\end{equation}
Note that the (strong) Markov property for the coupling implies the
(strong) marginal Markov property.

\section{Ergodicity of the environment process}\label{sectionEP}
\subsection{Assumptions on the environment}
In order to conclude results for the random walk, we need to have
sufficient control on the environment. To this end we assume there
exists a translation invariant coupling $\widehat{\bP}{}^E_{\eta,\xi}$ of the
environment, which satisfies the strong marginal Markov property (\ref
{eqSMMP}). In this coupling we look at $\widehat{\bE}{}^E_{\eta,\xi
}\rho(\eta
^1_t,\eta^2)$, measuring the distance of the states at the origin. If
this\vspace*{1pt} decays sufficiently fast we will be able to obtain ergodicity
properties of the environment.
\renewcommand{\theassumption}{1\textup{\alph{assumption}}}
\begin{assumption}\label{assumption1a}
The coupling $\widehat{\bP}{}^E$ satisfies
\[
\int_0^\infty t^d \sup
_{\eta,\xi\in\Omega} \widehat{\bE }{}^E_{\eta,\xi}\rho \bigl(
\eta_t^1(0),\eta_t^2(0)\bigr) \,dt<
\infty.
\]
\end{assumption}
This assumption is already sufficient to obtain the law of large
numbers for the position of the walker and unique ergodicity of the
environment process,
but it does not give quite enough control on local fluctuations. The
following stronger assumption remedies that.
\begin{assumption}\label{assumption1b}
The coupling $\widehat{\bP}{}^E$ satisfies
\[
\int_0^\infty t^d \sum
_{x\in\bZ^d}\sup_{(\eta,\xi)\in(\Omega
\times\Omega)_{0}} \widehat{
\bE}{}^E_{\eta,\xi}\rho\bigl(\eta _t^1(x),
\eta _t^2(x)\bigr) \,dt<\infty.
\]
\end{assumption}
\begin{remark*}
Typically, a coupling which satisfies Assumption~\ref{assumption1b}
also satisfies Assumption~\ref{assumption1a}. It is, however, not
automatic. But given a translation invariant coupling
$\widehat{\bP}{}^E $ which satisfies\vadjust{\goodbreak} Assumption~\ref{assumption1b} it
is possible to construct from $\widehat{\bP }{}^E$ a~new coupling
$\widetilde{\bP}{}^E$ via a telescoping argument so that
$\widetilde{\bP}$ satisfies both Assumptions~\ref{assumption1b} and
\ref{assumption1a}.
\end{remark*}

%
In Section~\ref{sectionexamples} we will discuss some examples which
satisfy those assumptions. Beside natural examples where $\widehat{\bE
}{}^E_{\eta,\xi}\rho(\eta_t^1(0),\eta_t^2(0))$ decays exponentially fast,
we give an example where other decay rates like polynomial decay are
obtained.

\subsection{Statement of the main theorem}
The main result of this section is the following theorem, which tells
us how the coupling property of the environment lifts to the
environment process.
%
\begin{theorem}\label{thmEP-estimate}
Let $f\dvtx \Omega\rightarrow\bR$ with $\interleave f \interleave<\infty$.
\begin{longlist}[(b)]
\item[(a)] Under Assumption~\ref{assumption1a}, there exists a constant $C_a>0$ so that
\[
\sup_{\eta,\xi\in\Omega}\int_0^\infty\bigl
\llvert S^{\mathrm{EP}}_t f(\eta)-S^{\mathrm{EP}}_t f(
\xi) \bigr\rrvert \,dt\leq C_a\interleave f \interleave.
\]
\item[(b)] Under Assumption~\ref{assumption1b}, there exists a constant $C_b>0$ so that
\[
\sum_{x\in\bZ^d}\sup_{(\eta,\xi)\in(\Omega\times\Omega)_{x}} \int
_0^\infty \bigl\llvert S^{\mathrm{EP}}_t
f(\eta)-S^{\mathrm{EP}}_t f(\xi) \bigr\rrvert \,dt\leq
C_b\interleave f \interleave.
\]
\end{longlist}
\end{theorem}
This theorem is the key to understanding the limiting behavior of the
random walk, that is, the law of large numbers, as well as for the
central limit theorem. Section~\ref{sectionproofs} is devoted the
proof of Theorem~\ref{thmEP-estimate}. In Section \ref
{sectionconvergence-speed} we generalize this result to give more
information about decay in time. Here we continue with results we can
obtain using Theorem~\ref{thmEP-estimate}. Most results about the
environment process just use part (a) of the theorem; part (b) shows
how more sophisticated properties lift from the environment to the
environment process as well. Those can be necessary to obtain more
precise results on the walker, like how likely atypical excursions from
the expected trajectory are.

It is possible to lift other properties from the environment to the
environment process as well. For example,
if Assumption~\ref{assumption1b} is modified to state
\[
\int_0^\infty t^d \sum
_{x\in\bZ^d}\sup_{(\eta,\xi)\in(\Omega
\times\Omega)_{0}} \widehat{
\bE}{}^E_{\eta,\xi}\frac{\rho(\eta
_t^1(x),\eta_t^2(x))}{\rho(\eta(0),\xi(0))} \,dt<\infty,
\]
then that implies for the environment process
\[
\sum_{x\in\bZ^d}\sup_{(\eta,\xi)\in(\Omega\times\Omega
)_{x}}\int
_0^\infty \frac{\llvert  S^{\mathrm{EP}}_t f(\eta)-S^{\mathrm{EP}}_t f(\xi) \rrvert
}{\rho(\eta(x),\xi(x))} \,dt\leq C_{b'}
\interleave f \interleave.
\]
This kind of condition can be relevant in the context of diffusive
environments to show that small changes
in the environment are causing only small changes in the environment process.

\subsection{Existence of a unique ergodic measure and continuity in
the rates}
First, the environment process, that is, the environment as seen from
the walker, is ergodic.
%
\begin{lemma}\label{lemmaergodicity}
Under Assumption~\ref{assumption1a} the environment process has a unique ergodic
probability measure $\mu^{\mathrm{EP}}$.
\end{lemma}
\begin{pf}
As $E$ is compact, so is $\Omega$, and therefore the space of
stationary measures is nonempty. So we must just prove uniqueness.

Assume $\mu,\nu$ are both stationary measures. Choose an arbitrary $f\dvtx
\Omega\rightarrow\bR$ with $\interleave f \interleave<\infty$. By
Theorem~\ref{thmEP-estimate}(a), for any $T>0$,
\begin{eqnarray*}
T\bigl\llvert \mu(f)-\nu(f) \bigr\rrvert &\leq& \iiint_0^T
\bigl\llvert S^{\mathrm{EP}}_t f (\eta)-S^{\mathrm{EP}}_tf(
\xi) \bigr\rrvert \,dt\, \mu (d\eta) \nu(d\xi)
\\
&\leq& \sup_{\eta,\xi\in\Omega} \int_0^\infty
\bigl\llvert S^{\mathrm{EP}}_t f (\eta)-S^{\mathrm{EP}}_t
f(\xi) \bigr\rrvert \,dt<\infty.
\end{eqnarray*}
As $T$ is arbitrary, $\mu(f)=\nu(f)$. As functions $f$ with
$\interleave f \interleave<\infty$ are dense in $\cC(\Omega)$,
there is at most one
stationary probability measure.
\end{pf}
It is of interest not only to know that the environment process has a
unique ergodic measure $\mu^{\mathrm{EP}}$, but also to know how this measure
depends on the rates $\alpha$.
%
\begin{theorem}\label{thmmu-continuous}
Under Assumption~\ref{assumption1a}, the unique ergodic measure $\mu^{\mathrm{EP}}_{\alpha}$
depends continuously on the rates $\alpha$. For two transition rate
functions $\alpha, \alpha'$, we have the following estimate:
\[
\bigl\llvert \mu^{\mathrm{EP}}_{\alpha}(f)-\mu^{\mathrm{EP}}_{\alpha'}(f)
\bigr\rrvert \leq\frac{C(\alpha)}{p(\alpha)}\bigl\llVert \alpha-\alpha' \bigr
\rrVert _0 \interleave f \interleave,
\]
that is,
\[
(\alpha,f) \mapsto\mu^{\mathrm{EP}}_\alpha(f)
\]
is continuous in $\llVert \cdot\rrVert _0\times\interleave
\cdot\interleave$.
The functions $C(\alpha), p(\alpha)$ satisfy $C(\alpha)>0, p(\alpha
)\in\ ]0,1[$. In the case that the rates $\alpha$ do not depend on
the environment, that is, $\alpha(\eta,z)=\alpha(z)$, they are given
by $p(\alpha)=1$,
\[
C(\alpha) = \int_0^\infty\sup
_{\eta,\xi\in\Omega} \widehat{\bE }{}^E_{\eta,\xi}\rho\bigl(
\eta^1_t(0),\eta^2_t(0)\bigr)\,dt.
\]
\end{theorem}
As the proof is a variation of the proof of Theorem \ref
{thmEP-estimate}, it is delayed to the end of Section~\ref{sectionproofs}.

\subsection{Speed of convergence to equilibrium in the environment
process}\label{sectionconvergence-speed}
We already know that under Assumption~\ref{assumption1a} the environment process has a
unique ergodic distribution. However, we do not know at what speed this
process converges to its unique stationary measure.\vadjust{\goodbreak} Given the speed of
convergence for the environment it is natural to believe that the
environment process inherits that speed with some form of slowdown due
to the additional self-interaction which is induced from the random
walk. For example, if the original speed of convergence were
exponential, then the environment process would also converge
exponentially fast. This is indeed the case.
%
\begin{theorem}\label{thmEP-phi-estimate}
Let $\phi\dvtx [0,\infty[\ \to\bR$ be a monotone increasing and
continuous function satisfying $\phi(0)=1$ and $\phi(s+t)\leq\phi
(s)\phi(t)$.
\begin{longlist}[(b)]
\item[(a)] Suppose the coupling $\widehat{\bP}{}^E$ satisfies
\[
\int_0^\infty\phi(t) t^d \sup
_{\eta,\xi\in\Omega} \widehat {\bE}{}^E_{\eta,\xi}\rho\bigl(
\eta_t^1(0),\eta_t^2(0)\bigr)\,dt<
\infty.
\]
Then there exists a constant $K_0>0$ and a decreasing function $C_a\dvtx
]K_0,\infty[\ \to[0,\infty[$ so that for any $K>K_0$ and any
$f\dvtx \Omega\rightarrow\bR$ with $\interleave f \interleave<\infty$,
\[
\sup_{\eta,\xi\in\Omega}\int_0^\infty\phi
\biggl(\frac
{t}{K} \biggr)\bigl\llvert S^{\mathrm{EP}}_t f(
\eta)-S^{\mathrm{EP}}_t f(\xi) \bigr\rrvert \,dt\leq
C_a(K)\interleave f \interleave.
\]
\item[(b)]Suppose the coupling $\widehat{\bP}{}^E$ satisfies
\[
\int_0^\infty\phi(t) t^d \sum
_{x\in\bZ^d}\sup_{(\eta,\xi)\in
(\Omega\times\Omega)_{0}} \widehat{
\bE}{}^E_{\eta,\xi}\rho\bigl(\eta _t^1(x),
\eta_t^2(x)\bigr)\,dt<\infty.
\]
Then there exists a constant $K_0>0$ and a decreasing function $C_b\dvtx
]K_0,\infty[\ \to[0,\infty[$ so that for any $K>K_0$ and any
$f\dvtx \Omega\rightarrow\bR$ with $\interleave f \interleave<\infty$,
\[
\sum_{x\in\bZ^d}\sup_{(\eta,\xi)\in(\Omega\times\Omega
)_{x}}\int
_0^\infty \phi \biggl(\frac{t}{K} \biggr)
\bigl\llvert S^{\mathrm{EP}}_t f(\eta)-S^{\mathrm{EP}}_t
f(\xi) \bigr\rrvert \, dt\leq C_b(K)\interleave f \interleave.
\]
\end{longlist}
\end{theorem}
Canonical choices for $\phi$ are $\phi(t)=\exp(\beta t^\alpha),
0<\alpha\leq1$ or $\phi(t)= (1+t)^\beta$, $\beta>0$.
This leads to the following transfer of convergence speed to
equilibrium from the environment to the environment process:
\begin{itemize}
\item exponential decay: $e^{-\lambda t} \longrightarrow
e^{-{\lambda t}/({K_0+\varepsilon}) }$,
\item stretched exponential decay: $e^{-\lambda t^\alpha}
\longrightarrow e^{-{\lambda t^\alpha}/{(K_0+\varepsilon)^\alpha}}$,
\item polynomial decay: $t^{-\lambda} \longrightarrow t^{-(\lambda
-d-\varepsilon)}$,
\end{itemize}
with $\varepsilon>0$ arbitrary, and in the case of polynomial decay,
$\lambda>d+1$.

\subsection{Consequences for the walker}
\label{sectionconsequences}
The strong convergence of the environment process to its stationary
measure obtained in Theorem~\ref{thmEP-estimate} implies various facts
for the random walker.
The most basic fact is that the random walker has a limiting speed.
%
\begin{proposition}\label{propLLN}
For any $\eta\in\Omega,x\in\bZ^d$,
\[
v:=\lim_{t\to\infty}\frac{X_t}{t} = \int\sum
_{z\in\bZ^d}z\alpha (\eta,z) \mu^{\mathrm{EP}}(d\eta)
\]
in $L^1$ and almost surely w.r.t. $\bP_{\eta,x}$. The
$L^1$-convergence is also uniform w.r.t. $\eta$ for a given $x$.
\end{proposition}
The convergence under $\bP_{\mu^{\mathrm{EP}},0}$ is a direct consequence of
ergodicity. For the extension to $\bP_{\eta,x}$ some ingredients of
the proofs in Section~\ref{sectionproofs} are required. Therefore the
proof is situated at the end of Section~\ref{sectionproofs}.

In the following theorem we prove the functional central limit theorem
for the position of the walker. The convergence to Brownian motion via
martingales is a rather straightforward consequence of the ergodicity
given by Theorem~\ref{thmEP-estimate}. The issue of nondegeneracy of
the variance is less standard and hence we give a proof.
%
\begin{theorem}
Assume Assumption~\ref{assumption1a}, $\llVert \alpha\rrVert _2<\infty$,
$\interleave\alpha\interleave_1<\infty$. Then the scaling limit of
the random walk is a Brownian
motion with drift $v$, that is,
\[
\frac{X_{tT}-vtT}{\sqrt{T}} \konv{T\to\infty} W_D(t),
\]
where $W_D$ is a Brownian motion with covariance matrix $D$.

Let $e\in\bR^d$ be a unit vector. Assume that either:
\begin{longlist}[(b)]
\item[(a)] there exists a $z\in\bZ^d, \langle e,z\rangle\neq0$, so
that for all $t>0$ and $\eta\in\Omega$ the probability $\bP_\eta
(\alpha(\eta_t,z)>0)$ is positive;
\item[(b)] $\mu^{\mathrm{EP}}(\alpha(\cdot,z))>0$ for $z\in\bZ^d$ with
$\langle e,z\rangle$ arbitrary large.
\end{longlist}
Then $\lim_{T\to\infty}\frac1T\Var(\langle X_T,e \rangle)>0$. In
particular, if \textup{(a)} or \textup{(b)} is satisfied for all $e$, then the covariance
matrix $D$ is nondegenerate.
\end{theorem}
\begin{pf}
Notice that $\sum_{z\in\bZ^d}z\alpha(\cdot,z)-v$ is in the domain
of $(L^{\mathrm{EP}})^{-1}$ because of Theorem~\ref{thmEP-estimate}. Decompose
\begin{eqnarray*}
X_t-vt &=& \biggl(X_t-\int_0^t
\sum_{z\in\bZ^d}z\alpha(\theta _{-X_s}
\eta_s,z)\,ds \biggr)
\\
&&{} + \biggl(\int_0^t \sum
_{z\in\bZ^d}z \bigl[\alpha(\theta _{-X_s}
\eta_s,z)-\mu^{\mathrm{EP}}\bigl(\alpha(\cdot,z)\bigr) \bigr]\,ds
\biggr).
\end{eqnarray*}
The first term on the right-hand side is a martingale, and the second
one is one as well, up to a uniformly bounded error. Both converge to
Brownian motion with finite variance by standard arguments when $\llVert \alpha\rrVert _2<\infty$. However, as the two terms are
not independent, an
argument is needed to prove that they do not annihilate. To prove that
we show that $\frac{1}{T}\Var(\langle X_T, e\rangle)$ is bounded
away from 0 under the assumed conditions.
Assume $T>0$ integer, and let $(\cF_t)_{t\geq0}$ be the canonical
filtration. Introduce the discrete-time martingale
\begin{eqnarray*}
M_n^T &:=& \bE\bigl[ \langle X_T, e\rangle|
\cF_n \bigr] - \bE\bigl[\langle X_T, e\rangle|
\cF_0\bigr] = X_n + \Psi_{T-n}(
\theta_{-X_n}\eta_n)-\Psi _t(
\eta_0),
\\
\Psi_S(\eta) &:=& \bE_{0,\eta}\int_0^{S}
\sum_{z\in\bZ^d}\langle z, e\rangle\alpha(
\theta_{-X_t}\eta_t,z) \,dt = \int_0^S
S_t^{\mathrm{EP}}\phi(\eta)\,dt;
\\
\phi(\eta)&:=& \sum_{z\in\bZ^d}\langle z, e\rangle\alpha(
\eta,z).
\end{eqnarray*}
With this, by stationarity of the environment process started from $\mu^{\mathrm{EP}}$,
\begin{eqnarray*}
\Var_{\mu^{\mathrm{EP}}}\bigl(\langle X_T, e\rangle\bigr) &\geq&
\bE_{\mu^{\mathrm{EP}}} \bigl(\langle X_T, e\rangle-\bE\bigl[\langle
X_T, e\rangle| \cF_0\bigr] \bigr)^2 \\
&=& \sum
_{n=1}^T \bE_{\mu^{\mathrm{EP}}}
(M_n-M_{n-1} )^2
\\
&=& \sum_{n=1}^T \bE_{\mu^{\mathrm{EP}}}
\bigl(\langle X_n, e\rangle-\langle X_{n-1}, e\rangle+
\Psi_{T-n}(\theta_{-X_n}\eta_n) \\
&&\hspace*{95.5pt}{} - \Psi
_{T-(n-1)}(\theta_{-X_{n-1}}\eta_{n-1}) \bigr)^2
\\
&=& \sum_{n=1}^T \bE_{\mu^{\mathrm{EP}}}
\bigl(\langle X_1, e\rangle+ \Psi _{T-n}(
\theta_{-X_1}\eta_1) - \Psi_{T-(n-1)}(
\eta_0) \bigr)^2.
\end{eqnarray*}
What has to be shown is that the above term is not 0. By Theorem \ref
{thmEP-estimate} and $\llVert \phi\rrVert _\infty\leq
\interleave\alpha\interleave_1<\infty$,
\[
\sup_{\eta,\xi\in\Omega} \sup_{T\geq0} \bigl\llvert
\Psi_{T}(\xi ) - \Psi_{T+1}(\eta) \bigr\rrvert =:C<\infty.
\]
Therefore, using $\llvert  a+b \rrvert \geq\llvert |a|-|b|
\rrvert $,
\[
\Var_{\mu^{\mathrm{EP}}}\bigl(\langle X_T, e\rangle\bigr) \geq T
\bE_{\mu^{\mathrm{EP}}}\ind _{\llvert \langle X_1, e\rangle\rrvert >C} \bigl(\bigl\llvert \langle
X_1, e\rangle\bigr\rrvert -C \bigr)^2.
\]
What remains to show is that $\bP_{\mu^{\mathrm{EP}}}(\llvert \langle X_1,
e\rangle\rrvert >C)>0$. If (b) is satisfied, this is immediate.
If (a) is
satisfied, then there is a positive probability that $X_t$ performs
sufficiently many jumps of size $z$ (and no other jumps) up to time~1.
\end{pf}
\begin{remark*}
The convergence to Brownian
motion with a nondegenerate variance also provides information about
the recurrence behavior of the walker. If $v=0$, supposing $d=1$ (in
higher dimensions, project onto a line), the limiting Brownian\vadjust{\goodbreak} motion
is centered. Hence there exists an infinite sequence $t_1<t_2<\cdots$
of times with $X_{t_{2n}}<0$ and $X_{t_{2n+1}}>0$, $n\in\bN$.
Supposing the walker has only jumps of size 1, it will traverse the
origin between $t_n,t_{n+1}$ for any $n\in\bN$; that is, it is
recurrent. (If the walker also has larger jumps, then one needs an
argument to actually hit the origin with some positive probability in
$[t_n,t_{n+1}]$.) Particularly, the recurrence implies that there
exists no regime where the random walk is transient but with 0 speed.
\end{remark*}

%
\section{Examples: Layered environments}
\label{sectionexamples}
There are many examples of environments which satisfy both
Assumptions~\ref{assumption1a} and~\ref{assumption1b}. Naturally, exponential convergence to the ergodic
measure is sufficient, independent of the dimension $d$. Therefore
interacting particle systems in the so-called $M<\varepsilon$-regime or
weakly interacting diffusions on a compact manifold belong to the
environments to which this method is applicable.

To exploit the fact that only sufficient polynomial decay of
correlations is required, we will construct a class of environments
which we call \textit{layered environments}. One can think of layered
environments as a weighted superposition of a sequence of (independent)
environments.

Those kind of environments are fairly natural objects to study. One
area where they can appear is an idealization of molecular motors. In
molecular motors the walker moves (e.g.) in a potential, where the
potential randomly switches between various global states [e.g.,
related to chemical transitions in the example of kinesine; see \citet
{AMBAYEKEHR99,JARZYNSKIMAZONKA99,JULICHERAJDARIPROST97,MAGNASCO94,DONATOPIATNITSKI05}
for more motivation]. Here each layer is representing the interaction
with the environment for one global state. In many realistic situations
there are many such states. If the global state changes very quickly
compared to the movement of the random walk, what is observed is a
weighted superposition with weights given by the relative frequencies
of the appearance of the individual global states.

Layered environments could also appear from a multi-scale analysis of a
complicated environment, where the layers with a high index represent
the long-range interactions. Besides, layered environments are a useful
tool because they form a class of environments which are uniformly
mixing with arbitrary mixing speed. There are plenty of examples where
one has polynomial or stretched exponential mixing, for example, in the
context of diffusion processes. However, those examples are not
uniformly mixing, in the context of diffusion processes because of an
unbounded state space.

Here we focus on layers which still have exponential decay of
correlations, but each layer does converge to its stationary measure at
a layer specific rate $\alpha_n$, with $n$ being the index of the
layer. When $\alpha_n$ tends to 0 as $n\to\infty$ this introduces
some form of arbitrary slow decay of correlations. We counterbalance
this by weighting the superposition in such a way that the individual
influence of a layer goes to 0 as well. Note that such a
counterbalancing is only possible\vadjust{\goodbreak} because of the Lipschitz nature of
the assumptions. A uniform decay estimate does not hold because of the
arbitrary slow decay in deep layers.

More formally, for each $n\in\bN$ let $(\eta^n)_{t\geq0}$ be a
Markov process on $\Omega_0:=\{0,1\}^{\bZ^d}$, the environment on
layer $n$. This process should have a coupling $\widehat{\bP}{}^n_{\eta,\xi}$ with
%
\begin{equation}
\label{eqlayerdecay} \sup_{\eta,\xi\in\Omega_0}\widehat{\bE}{}^n_{\eta,\xi}
\bigl\llvert \eta ^{n,1}_t(0)-\eta^{n,2}_t(0)
\bigr\rrvert \leq2 e^{-\alpha_n t},\qquad \alpha_n>0.
\end{equation}
The \textit{layered environment} $(\eta_t)_{t\geq0}$ then consists of
the stack of independent layers $(\eta^n_t)_{t\geq0}$. The single
site state space is $E=\{0,1\}^\bN$ and space of all configurations
$\Omega=E^{\bZ^d}$.

The superposition of the environments is weighted by the distance $\rho
$ on $E$, which we choose in the following way. Fix a sequence $\gamma
_1>\gamma_2>\cdots>0$ with $\sum_{n\in\bN}\gamma_n=1$. For
$(a_n)_{n\in\bN},(b_n)_{n\in\bN}\in E$ the distance is
%
\begin{equation}
\label{eqstackdistance} \rho\bigl((a_n), (b_n)\bigr):=
\sum_{n\in\bN} \gamma_n\llvert
a_n-b_n \rrvert.
\end{equation}
The coupling $\widehat{\bP}{}^E$ of the layered environments is simply the
independent coupling of the individual layer couplings $\widehat{\bP
}{}^n$. The
layer decay (\ref{eqlayerdecay}) and the choice of distance (\ref
{eqstackdistance}) then provide the following decay of coupling distance
for the layered environment:
%
\begin{equation}
\label{eqstackdecay} \sup_{\eta,\xi\in\Omega}\widehat{\bE}{}^E_{\eta,\xi}
\rho\bigl(\eta ^{1}_t(0),\eta^{2}_t(0)
\bigr) \leq2 \sum_{n\in\bN}\gamma_n
e^{-\alpha
_n t}.
\end{equation}
The sum on the right-hand side of (\ref{eqstackdecay}) can have
arbitrary slow decay depending on $\alpha_n,\gamma_n$. For example,
if one fixes $\alpha_n=n^{-1}$, then $\gamma_n = n^{-\gamma-1}$
leads to decay of order $t^{-\gamma}$.

We did not specify the exact nature of the individual layers, as those
did not matter for the construction. A natural example is when
individual layers consist of Ising model Glauber dynamics at inverse
temperature $\beta_n<\beta_c$, and $\beta_n \to\beta_c$ as $n\to
\infty$.

\section{Proofs}\label{sectionproofs}
In this section we always assume that Assumption~\ref{assumption1a} holds.

We start with an outline of the idea of the proofs. We have a coupling
of the environments $(\eta^1_t,\eta^2_t)$, which we extend to include
two random walkers $(X^1_t,X^2_t)$, driven by their corresponding
environment. We maximize the probability of both walkers performing the
same jumps. Then Assumption~\ref{assumption1a} is sufficient to obtain a positive
probability of both walkers staying together forever. If the walkers
stay together, one just has to account for the difference in
environments, but not the walkers as well. When the walkers split, the
translation invariance allows for everything to shift so that both
walkers are back at the origin, and one can try again. After a
geometric number of trials it is then guaranteed that the walkers stay together.
%
\begin{proposition}[(Coupling construction)]\label{propcoupling-construction}
Given the coupling $\widehat{\bP}{}^E_{\eta,\xi}$ of the
environments, we extend
it to a coupling $\widehat{\bP}_{\eta,x;\xi,y}$. This coupling has the
following properties:
\begin{enumerate}[(b)]
\item[(a)] (Marginals) The coupling supports two environments and
corresponding random walkers:
\begin{enumerate}[(2)]
\item[(1)] $\widehat{\bP}_{\eta,x;\xi,y}((\eta^1_t,X_t^1)\in\cdot) = \bP
_{\eta,x}((\eta_t,X_t) \in\cdot)$;
\item[(2)] $\widehat{\bP}_{\eta,x;\xi,y}((\eta^2_t,X_t^2)\in\cdot) = \bP
_{\xi,y}((\eta_t,X_t) \in\cdot)$;
\end{enumerate}
\item[(b)] (Extension of $\widehat{\bP}{}^E_{\eta,\xi}$) The
environments behave
as under $\widehat{\bP}{}^E$,
\[
\widehat{\bP}_{\eta,x;\xi,y}\bigl(\bigl(\eta^1_t,
\eta_t^2\bigr)\in\cdot\bigr) = \widehat{\bP
}{}^E_{\eta,\xi}\bigl(\bigl(\eta_t^1,
\eta_t^2\bigr) \in\cdot\bigr);
\]
\item[(c)]\label{enumdecoupling-rate} (Coupling of the walkers)
$X_t^1$ and $X_t^2$ perform identical jumps as much as possible, and
the rate of performing a different jump is $\sum_{z\in\bZ^d}\llvert \alpha(\theta_{-X_t^1}\eta_t^1,\break z)-\alpha(\theta_{-X_t^2}\eta
_t^2,z) \rrvert $;
\item[(d)] (Minimal and maximal walkers) In addition to the
environments $\eta_t^1$ and $\eta_t^2$ and random walkers $X_t^1$ and
$X_t^2$, the coupling supports minimal and maximal walkers $Y_t^+$,
$Y_t^-$ as well. These two walkers have the following properties:
\begin{enumerate}[(4)]
\item[(1)] $Y_t^- \leq X_t^1-x, X_t^2-y \leq Y_t^+
\widehat{\bP}_{\eta,x;\xi,y}$-a.s. (in dimension $d>1$, this is to be
interpreted coordinate-wise);\vspace*{1pt}
\item[(2)] $Y_t^+,Y_t^-$ are independent of
$\eta_t^1,\eta_t^2$;\vspace*{1pt}
\item[(3)] $\widehat{\bE}_{\eta,x;\xi,y}Y_t^+ = t\gamma^+$ for some
$\gamma
^+\in\bR^d$;\vspace*{1pt}
\item[(4)] $\widehat{\bE}_{\eta,x;\xi,y}Y_t^- = t\gamma^-$ for some
$\gamma
^-\in\bR^d$.
\end{enumerate}
\end{enumerate}
\end{proposition}
\begin{pf}
The construction of this coupling $\widehat{\bP}_{\eta,x;\xi,y}$ is
done in
the following way: we extend the original coupling $\widehat{\bP
}{}^E_{\eta,\xi}$
to contain an independent sequence of Poisson processes $N^z,z\in\bZ
^d$, with rates $\lambda_z:=\sup_\eta\alpha(\eta,z)$, as well as a
sufficient supply of independent uniform $[0,1]$ variables. The walkers
$X^1, X^2$ then start from~$x$ (resp., $y$) and exclusively (but not
necessarily) jump when one of the Poisson clocks $N^z$ rings. When the
clock $N^z$ rings the walkers jumps from $X^i_t$ to $X^i_t+z$ only if a
uniform $[0,1]$ variable $U$ satisfies $U<\alpha(\theta_{-X_t^i}\eta
_t^i,z)/\lambda_z$, $i=1,2$. Note that both walkers share the same
$U$, but $U$'s for different rings of the Poisson clocks are independent.

The upper and lower walkers $Y_t^+,Y_t^-$ are constructed from the same
Poisson clocks $N^z$.
They always jump on these clocks; however, they jump by $\max(z,0)$ or
$\min(z,0)$, respectively.

The properties of the coupling arise directly from the construction
plus the fact that $\llVert \alpha\rrVert _1<\infty$.
\end{pf}

To ease notation we will call $\widehat{\bP}_{\eta,0;\xi,0}$ simply
$\widehat{\bP}
_{\eta,\xi}$ and the law of $Y_t^+,Y_t^-$ $\widehat{\bP}$ whenever
there is no
fear of confusion.\vadjust{\goodbreak}

Now we show how suitable estimates on the coupling speed of the
environment translate to properties of the extended coupling.\vspace*{-2pt}
%
\begin{lemma}\label{lemmasupnorm-coupling-estimate}
\[
\widehat{\bE}_{\eta,x;\xi,y}\rho\bigl(\eta_t^1
\bigl(X_t^1\bigr), \eta _t^2
\bigl(X_t^1\bigr)\bigr)<\bigl(\bigl\llVert
\gamma^{+} - \gamma^{-} \bigr\rrVert _\infty t + 1
\bigr)^d\sup_{\eta,\xi\in\Omega
}\widehat{\bE}{}^E_{\eta,\xi}
\rho\bigl(\eta_t^1(0),\eta_t^2(0)
\bigr).\vspace*{-2pt}
\]
\end{lemma}
\begin{pf}
Denote with $R_t \subset\bZ^d$ the set of sites $z\in\bZ^d$ with
$Y_t^-\leq z \leq Y_t^+$ (coordinate-wise). Then
\begin{eqnarray*}
&&
\sup_{\eta,\xi,x,y} \widehat{\bE}_{\eta,x;\xi,y}\rho\bigl(
\eta_t^1\bigl(X_t^1\bigr),
\eta_t^2\bigl(X_t^1\bigr)\bigr)
\\[-2pt]
&&\qquad\leq \sup_{\eta,\xi,x,y} \widehat{\bE}_{\eta,x;\xi,y}\sum
_{z\in
R_t}\rho\bigl(\eta_t^1(x+z),
\eta_t^2(x+z)\bigr)
\\[-2pt]
&&\qquad\leq \widehat{\bE} \biggl[\sum_{z\in R_t}1 \biggr]\sup
_{\eta,\xi,z} \widehat{\bE}{}^E_{\eta,\xi}\rho\bigl(
\eta_t^1(z), \eta_t^2(z)\bigr)
\\[-2pt]
&&\qquad \leq\bigl(\bigl\llVert \gamma^{+} - \gamma^{-} \bigr
\rrVert _\infty t + 1\bigr)^d \sup_{\eta,\xi\in\Omega}
\widehat{\bE}{}^E_{\eta,\xi}\rho\bigl(\eta _t^1(0),
\eta_t^2(0)\bigr).\vspace*{-2pt}
\end{eqnarray*}
\upqed\end{pf}
%
\begin{lemma}\label{lemmadecoupling-estimate}
Denote by $\tau:=\inf\{t\geq0\dvtx X^1_t\neq X_t^2\}$ the first time the
two walkers are not at the same position. Under Assumption~\ref{assumption1a},
\[
\inf_{\eta,\xi\in\Omega} \widehat{\bP}_{\eta;\xi}(\tau=\infty) >0,
\]
that is, the walkers $X^1$ and $X^2$ never decouple with strictly
positive probability.\vspace*{-2pt}
\end{lemma}
\begin{pf}
Both walkers start in the origin, therefore $\tau>0$. The probability
that a Poisson clock with time dependent rate $\lambda_t$ is has not
yet rung by time $T$ is $\exp({-\int_0^T\lambda_t\,dt})$. As the
rate of decoupling is given by Proposition \ref
{propcoupling-construction}(c), we obtain
%
\begin{eqnarray}\label{eqdynamic-rate-poisson}\quad
\widehat{\bP}_{\eta,\xi}(\tau> T) &=& \widehat{\bE }_{\eta,\xi} \exp
\biggl({-\int_0^T \sum
_{z\in\bZ^d}\bigl\llvert \alpha\bigl(\theta _{-X^1_t}
\eta^1_t,z\bigr)- \alpha\bigl(\theta_{-X^1_t}
\eta^2_t,z\bigr) \bigr\rrvert \,dt} \biggr)
\nonumber\\[-9pt]\\[-9pt]
&\geq&\exp \biggl(-\widehat{\bE}_{\eta,\xi}{\int_0^T
\sum_{z\in\bZ
^d}\bigl\llvert \alpha\bigl(
\theta_{-X^1_t}\eta^1_t,z\bigr)- \alpha\bigl(\theta
_{-X^1_t}\eta^2_t,z\bigr) \bigr\rrvert \,dt}
\biggr).
\nonumber
\end{eqnarray}
By telescoping over single site changes,
\begin{eqnarray*}
&&
\widehat{\bE}_{\eta,\xi}\sum_{z\in\bZ^d}\bigl\llvert
\alpha\bigl(\theta _{-X^1_t}\eta^1_t,z\bigr)-
\alpha\bigl(\theta_{-X^1_t}\eta^2_t,z\bigr) \bigr
\rrvert
\\[-2pt]
&&\qquad \leq\widehat{\bE}_{\eta,\xi}\sum_{z\in\bZ^d}\sum
_{x\in\bZ
^d}\rho\bigl(\eta^1_t
\bigl(X^1_t+x\bigr),\eta^2_t
\bigl(X^1_t+x\bigr)\bigr)\delta_{\alpha(\cdot,z)}(x)
\\
&&\qquad \leq\sup_{x\in\bZ^d}\widehat{\bE}_{\eta,\xi}\rho\bigl(\eta
^1_t\bigl(X^1_t+x\bigr),
\eta^2_t\bigl(X^1_t+x\bigr)\bigr)
\interleave\alpha\interleave
\\
&&\qquad \leq\interleave\alpha\interleave\bigl(\bigl\llVert \gamma^{+} -
\gamma^{-} \bigr\rrVert _\infty t + 1\bigr)^d\sup
_{\eta,\xi\in\Omega}\widehat{\bE}{}^E_{\eta,\xi}\rho \bigl(\eta
_t^1(0),\eta_t^2(0)\bigr),
\end{eqnarray*}
where the last line follows from Lemma \ref
{lemmasupnorm-coupling-estimate}. With this estimate and Assumption~\ref{assumption1a}
we obtain
\begin{eqnarray*}
&&
\widehat{\bP}_{\eta,\xi}(\tau=\infty)
\\
&&\qquad\geq\exp \biggl(- \interleave\alpha\interleave\int_0^\infty
\bigl(\bigl\llVert \gamma^{+} - \gamma^{-} \bigr\rrVert
_\infty t + 1\bigr)^d\sup_{\eta,\xi\in\Omega}\widehat{
\bE}{}^E_{\eta,\xi}\rho\bigl(\eta _t^1(0),
\eta_t^2(0)\bigr)\,dt \biggr)
\\
&&\qquad>0 \qquad\mbox{uniformly in $\eta, \xi$}.
\end{eqnarray*}
\upqed\end{pf}
\begin{pf*}{Proof of Theorem~\ref{thmEP-estimate}, part \textup{(a)}}
The idea of the proof is to use the coupling of Proposition \ref
{propcoupling-construction}: we wait until the walkers $X_t^1$ and
$X_t^2$, which are initially at the same position, decouple, and then
restart everything and try again. By Lem\-ma~\ref
{lemmadecoupling-estimate} there is a positive probability of never
decoupling, so this scheme is successful. Using the time of decoupling
$\tau$ (as in Lemma~\ref{lemmadecoupling-estimate}) and the strong
marginal Markov property (\ref{eqSMMP}),
%
\begin{eqnarray}
&&
\int_0^T \bigl\llvert \widehat{
\bE}_{\eta,0;\xi,0}\ind_{t\geq\tau} \bigl(f\bigl(\theta_{-X^1_t}
\eta^1_t\bigr)-f\bigl(\theta_{-X^2_t}
\eta^2_t\bigr) \bigr) \bigr\rrvert \, dt
\nonumber
\\
&&\qquad = \int_0^T \bigl\llvert \widehat{
\bE}_{\eta0;\xi,0}\ind_{t\geq\tau
}\erwsymbol \bigl[f\bigl(
\theta_{-X^1_t}\eta^1_t\bigr)-f\bigl(
\theta_{-X^2_t}\eta ^2_t\bigr)|\fF_\tau
\bigr] \bigr\rrvert \,dt
\nonumber
\\
&&\qquad \leq\int_0^T \widehat{\bE}_{\eta,0;\xi,0}
\ind_{t\geq\tau}\bigl\llvert S^{\mathrm{EP}}_{t-\tau}f\bigl(
\theta_{-X^1_\tau}\eta^1_\tau \bigr)-S^{\mathrm{EP}}_{t-\tau}f
\bigl(\theta_{-X^2_\tau}\eta^2_\tau\bigr) \bigr\rrvert
\, dt
\nonumber
\\
\label{eqphi-1}
&&\qquad = \widehat{\bE}_{\eta,0;\xi,0} \int_0^{(T-\tau)\vee0}
\bigl\llvert S^{\mathrm{EP}}_{t}f\bigl(\theta_{-X^1_\tau}
\eta^1_\tau\bigr)-S^{\mathrm{EP}}_{t}f\bigl(
\theta _{-X^2_\tau}\eta^2_\tau\bigr) \bigr\rrvert \,dt
\\
\label{eqphi-2}
&&\qquad \leq\widehat{\bP}_{\eta,0;\xi,0}(\tau<\infty) \sup_{\eta,\xi\in
\Omega}
\int_0^{T} \bigl\llvert S^{\mathrm{EP}}_{t}f(
\eta)-S^{\mathrm{EP}}_{t}f(\xi) \bigr\rrvert \,dt.
\end{eqnarray}
And therefore
%
\begin{eqnarray}
&&
\int_0^T \bigl\llvert S^{\mathrm{EP}}_t
f(\eta)-S^{\mathrm{EP}}_t f(\xi) \bigr\rrvert \,dt
\nonumber
\\
&&\qquad = \int_0^T \bigl\llvert \widehat{
\bE}_{\eta,\xi}f\bigl(\theta_{-X^1_t}\eta ^1_t
\bigr)-f\bigl(\theta_{-X^2_t}\eta^2_t\bigr) \bigr
\rrvert \,dt
\nonumber
\\
&&\qquad \leq\int_0^T \widehat{\bE}_{\eta,\xi}
\ind_{t<\tau}\bigl\llvert f\bigl(\theta_{-X^1_t}
\eta^1_t\bigr)-f\bigl(\theta_{-X^1_t}
\eta^2_t\bigr) \bigr\rrvert \,dt
\nonumber
\\
\label{eqphi-3}
&&\qquad\quad{} + \widehat{\bP}_{\eta,\xi}(\tau<\infty) \sup_{\eta,\xi\in\Omega
}\int
_0^{T} \bigl\llvert S^{\mathrm{EP}}_{t}f(
\eta)-S^{\mathrm{EP}}_{t}f(\xi) \bigr\rrvert \,dt
\\
&&\qquad \leq\int_0^\infty\widehat{\bE}_{\eta,\xi}
\bigl\llvert f\bigl(\theta _{-X^1_t}\eta^1_t
\bigr)-f\bigl(\theta_{-X^1_t}\eta^2_t\bigr) \bigr
\rrvert \,dt
\nonumber
\\
\label{eqphi-4}
&&\qquad\quad{} + \widehat{\bP}_{\eta,\xi}(\tau<\infty)\sup_{\eta,\xi\in\Omega
}\int
_0^T \bigl\llvert S^{\mathrm{EP}}_{t}f(
\eta)-S^{\mathrm{EP}}_{t}f(\xi) \bigr\rrvert \,dt,
\end{eqnarray}
which gives us the upper bound
%
\begin{eqnarray}\label{eqphi-5}
&&
\sup_{\eta,\xi\in\Omega}\int_0^\infty\bigl
\llvert S^{\mathrm{EP}}_tf(\eta)-S^{\mathrm{EP}}_tf(\xi)
\bigr\rrvert \,dt
\nonumber\\
&&\qquad \leq \Bigl(\inf_{\eta,\xi\in\Omega}\widehat{\bP}_{\eta,\xi}(\tau =
\infty) \Bigr)^{-1}\\
&&\qquad\quad{}\times\sup_{\eta,\xi\in\Omega}\int_0^\infty
\widehat{\bE} _{\eta,\xi}\bigl\llvert f\bigl(\theta_{-X^1_t}
\eta^1_t\bigr)-f\bigl(\theta _{-X^1_t}
\eta^2_t\bigr) \bigr\rrvert \,dt.\nonumber
\end{eqnarray}
To show that the last integral is finite, we telescope over single site
changes, and get
\begin{eqnarray*}
&&\int_0^\infty\widehat{\bE}_{\eta,\xi}\bigl
\llvert f\bigl(\theta_{-X^1_t}\eta ^1_t\bigr)-f
\bigl(\theta_{-X^1_t}\eta^2_t\bigr) \bigr\rrvert
\,dt
\\
&&\qquad \leq\int_0^\infty\widehat{\bE}_{\eta,\xi}
\sum_{x\in\bZ^d}\rho \bigl(\eta^1_t
\bigl(x+X^1_t\bigr),\eta^2_t
\bigl(x+X^1_t\bigr)\bigr)\delta_f(x)\,dt
\\
&&\qquad \leq\interleave f \interleave\sup_{\eta,\xi,x}\int
_0^\infty\widehat{\bE}_{\eta,\xi}\rho\bigl(
\eta^1_t\bigl(x+X_t^1\bigr),\eta
^2_t\bigl(x+X_t^1\bigr)\bigr)
\,dt,
\end{eqnarray*}
which is finite by Lemma~\ref{lemmasupnorm-coupling-estimate} and
Assumption~\ref{assumption1a}. Choosing
%
\begin{equation}\label{eqphi-6}\qquad
C_a = \Bigl(\inf_{\eta,\xi\in\Omega}\widehat{
\bP}_{\eta,\xi}(\tau =\infty) \Bigr)^{-1}\sup_{\eta,\xi,x}
\int_0^\infty\widehat{\bE}_{\eta,\xi} \rho
\bigl(\eta^1_t\bigl(x+X_t^1\bigr),
\eta^2_t\bigl(x+X_t^1\bigr)\bigr)
\,dt
\end{equation}
completes the proof.
\end{pf*}
To prove part (b) of the theorem, we need the following analogue to
Lemma~\ref{lemmadecoupling-estimate} using Assumption~\ref{assumption1b}.
%
\begin{lemma}\label{lemmatriplenorm-coupling-estimate}
Under Assumption~\ref{assumption1b}, for every site-weight function $w\dvtx\break \bZ^d \to
[0,\infty[$ with $\llVert  w \rrVert _1:= \sum_x
w(x)<\infty$,
we have
\begin{eqnarray*}
&&
\sum_{x\in\bZ^d}\sup_{(\eta,\xi)\in(\Omega\times\Omega
)_{x}}\int
_0^\infty \sum_{y\in\bZ^d}
w(y)\widehat{\bE}_{\eta,\xi}\rho \bigl(\eta_t^1
\bigl(y+X_t^1\bigr), \eta_t^2
\bigl(y+X_t^1\bigr)\bigr)\,dt\\
&&\qquad\leq \mathrm{const} \cdot\llVert w
\rrVert _1.
\end{eqnarray*}
\end{lemma}
\begin{pf}
Denote with $R_t \subset\bZ^d$ the set of sites whose $j$th
coordinate lies between $Y_t^{j,-}$ and $Y_t^{j,+}$. Then
\begin{eqnarray*}
&& \sum_{y\in\bZ^d} w(y)\widehat{\bE}_{\eta,\xi}\rho
\bigl(\eta_t^1\bigl(y+X_t^1\bigr),
\eta_t^2\bigl(y+X_t^1\bigr)\bigr)
\\
&&\qquad = \sum_{y\in\bZ^d} \widehat{\bE}_{\eta,\xi}w
\bigl(y-X_t^1\bigr)\rho\bigl(\eta _t^1(y),
\eta_t^2(y)\bigr)
\\
&&\qquad \leq\sum_{y\in\bZ^d} \widehat{\bE}_{\eta,\xi}\sum
_{z\in
R_t}w(y-z) \rho\bigl(\eta_t^1(y),
\eta_t^2(y)\bigr)
\\
&&\qquad = \sum_{y\in\bZ^d} \widehat{\bE} \biggl[\sum
_{z\in R_t}w(y-z) \biggr]\widehat{\bE}{}^E_{\eta,\xi}
\rho\bigl(\eta_t^1(y), \eta_t^2(y)
\bigr)
\end{eqnarray*}
by independence of $R_t$ and $(\eta^1_t,\eta^2_t)$. Therewith,
\begin{eqnarray*}
&&\sum_{x\in\bZ^d}\sup_{(\eta,\xi)\in(\Omega\times\Omega
)_{x}}\int
_0^\infty \sum_{y\in\bZ^d}
w(y)\widehat{\bE}_{\eta,\xi}\rho \bigl(\eta_t^1
\bigl(y+X_t^1\bigr), \eta_t^2
\bigl(y+X_t^1\bigr)\bigr)\,dt
\\
&&\qquad \leq\int_0^\infty\sum
_{y\in\bZ^d} \widehat{\bE} \biggl[\sum_{z\in
R_t}w(y-z)
\biggr] \sum_{x\in\bZ^d}\sup_{(\eta,\xi)\in(\Omega
\times\Omega)_{x}}
\widehat{\bE}{}^E_{\eta,\xi} \rho\bigl(\eta_t^1(y),
\eta _t^2(y)\bigr)\,dt.
\end{eqnarray*}
Note that by translation invariance the right-hand side is equal to
\[
\sum_{x\in\bZ^d}\sup_{(\eta,\xi)\in(\Omega\times\Omega
)_{0}} \widehat{
\bE}{}^E_{\eta,\xi} \rho\bigl(\eta_t^1(x),
\eta_t^2(x)\bigr).
\]
By construction of $R_t$ and Proposition \ref
{propcoupling-construction}(d),
\begin{eqnarray*}
\sum_{y\in\bZ^d} \widehat{\bE} \biggl[\sum
_{z\in R_t}w(y-z) \biggr] &=& \widehat{\bE} \biggl[\sum
_{z\in R_t} 1 \biggr]\llVert w \rrVert _1 = \prod
_{j=1}^d\bigl(\gamma^{j,+}t -
\gamma^{j,-}t+1\bigr) \llVert w \rrVert _1
\\
&\leq& c \bigl(t^d+1\bigr)\llVert w \rrVert _1
\end{eqnarray*}
for some suitable $c>0$.
Therefore Assumption~\ref{assumption1b} completes the proof.
\end{pf}
\begin{pf*}{Proof of Theorem~\ref{thmEP-estimate}, part \textup{(b)}}
Let $\tau:=\inf\{t\geq0\dvtx X^1_t\neq X^2_t\}$. Then we split the
integration at $\tau$,
\begin{eqnarray*}
&&\sum_{x\in\bZ^d}\sup_{(\eta,\xi)\in(\Omega\times\Omega
)_{x}}\int
_0^\infty \bigl\llvert S_t^{\mathrm{EP}}
f(\eta)-S_t^{\mathrm{EP}} f(\xi) \bigr\rrvert \,dt
\\
&&\qquad \leq\sum_{x\in\bZ^d}\sup_{(\eta,\xi)\in(\Omega\times
\Omega)_{x}}\int
_0^\infty\bigl\llvert \widehat{\bE}_{\eta,\xi}
\ind_{\tau
>t} \bigl(f\bigl(\theta_{-X_t^1}\eta_t^1
\bigr)-f\bigl(\theta_{-X_t^1}\eta_t^2\bigr) \bigr)
\bigr\rrvert \,dt
\\
&&\qquad\quad{} +\sum_{x\in\bZ^d}\sup_{(\eta,\xi)\in(\Omega\times
\Omega)_{x}}\int
_0^\infty\bigl\llvert \widehat{\bE}_{\eta,\xi}
\ind_{\tau\leq t} \bigl(f\bigl(\theta_{-X_t^1}\eta_t^1
\bigr)-f\bigl(\theta_{-X_t^2}\eta_t^2\bigr) \bigr)
\bigr\rrvert \,dt.
\end{eqnarray*}
We estimate the first term by moving the expectation out of the
absolute value and forgetting the restriction to $\tau>t$,
\[
\sum_{x\in\bZ^d}\sup_{(\eta,\xi)\in(\Omega\times\Omega
)_{x}}\int
_0^\infty \sum_{y\in\bZ^d}
\delta_f(y)\widehat{\bE}_{\eta,\xi
}\rho\bigl(
\eta_t^1\bigl(y+X_t^1\bigr),
\eta_t^2\bigl(y+X_t^1\bigr)\bigr)
\,dt.
\]
By Lemma~\ref{lemmatriplenorm-coupling-estimate} with $w = \delta_f$,
this is bounded by some constant times $\interleave f \interleave$.
For the second term we start by using the strong marginal Markov
property~(\ref{eqSMMP}),
%
\begin{eqnarray}
&& \int_0^\infty\bigl\llvert \widehat{
\bE}_{\eta,\xi}\ind_{\tau\leq t} \bigl(f\bigl(\theta_{-X_t^1}
\eta_t^1\bigr)-f\bigl(\theta_{-X_t^2}
\eta_t^2\bigr) \bigr) \bigr\rrvert \,dt
\nonumber
\\[-2pt]
&&\qquad = \int_0^\infty\bigl\llvert \widehat{
\bE}_{\eta,\xi}\ind_{\tau\leq
t} \bigl(S_{t-\tau}^{\mathrm{EP}}f
\bigl(\theta_{-X_\tau^1}\eta_\tau ^1
\bigr)-S_{t-\tau}^{\mathrm{EP}}f\bigl(\theta_{-X_\tau^2}
\eta_\tau^2\bigr) \bigr) \bigr\rrvert \,dt
\nonumber
\\[-2pt]
\label{eqphi-7}
&&\qquad \leq\widehat{\bE}_{\eta,\xi}\ind_{\tau<\infty}\int_\tau^\infty
\bigl\llvert \bigl(S_{t-\tau}^{\mathrm{EP}}f\bigl(\theta_{-X_\tau^1}
\eta_\tau ^1\bigr)-S_{t-\tau}^{\mathrm{EP}}f\bigl(
\theta_{-X_\tau^2}\eta_\tau^2\bigr) \bigr) \bigr\rrvert
\,dt
\\[-2pt]
\label{eqphi-8}
&&\qquad \leq \widehat{\bP}_{\eta,\xi}(\tau<\infty)\sup_{\eta,\xi\in\Omega}\int
_0^\infty\bigl\llvert S_t^{\mathrm{EP}}
f(\eta)-S_t^{\mathrm{EP}} f(\xi) \bigr\rrvert \, dt.
\end{eqnarray}
By part (a) of Theorem~\ref{thmEP-estimate} the integral part is
uniformly bounded by $C_a \interleave f \interleave$. So what remains
to complete the
proof is to show that
%
\begin{equation}\label{eqphi-9}
\sum_{x\in\bZ^d}\sup_{(\eta,\xi)\in(\Omega\times\Omega
)_{x}}\widehat{
\bP}_{\eta,\xi
}(\tau<\infty)< \infty.
\end{equation}

To do so we first use the same idea as in the proof of Lemma \ref
{lemmadecoupling-estimate} to obtain
\begin{eqnarray*}
&&
\widehat{\bP}_{\eta,\xi}(\tau<\infty)
\\[-2pt]
&&\qquad = 1-\exp \biggl(- \int_0^\infty\widehat{
\bE}_{\eta,\xi} \sum_{z\in
\bZ^d}\bigl\llvert \alpha
\bigl(\theta_{-X_t^1}\eta_t^1,z\bigr) - \alpha
\bigl(\theta_{-X_t^1}\eta_t^2,z\bigr) \bigr\rrvert
\,dt \biggr)
\\[-2pt]
&&\qquad \leq\int_0^\infty\widehat{\bE}_{\eta,\xi}
\sum_{z\in\bZ^d}\bigl\llvert \alpha\bigl(
\theta_{-X_t^1}\eta_t^1,z\bigr) - \alpha\bigl(
\theta _{-X_t^1}\eta_t^2,z\bigr) \bigr\rrvert \,dt
\\[-2pt]
&&\qquad \leq\int_0^\infty\sum
_{y\in\bZ^d} w_\alpha(y) \widehat{\bE} _{\eta,\xi}\rho
\bigl(\eta_t^1\bigl(y+X_t^1\bigr),
\eta_t^2\bigl(y+X_t^1\bigr)\bigr)
\,dt
\end{eqnarray*}
with
\[
w_\alpha(x):= \sup_{(\eta,\xi)\in(\Omega\times\Omega)_{x}} \sum
_{z\in\bZ
^d}\bigl\llvert \alpha(\eta,z) - \alpha(\xi,z) \bigr\rrvert
\]
and $\sum_{x\in\bZ^d}w_\alpha(x)<\infty$.
So we get
\begin{eqnarray*}
&&
\sum_{x\in\bZ^d}\sup_{(\eta,\xi)\in(\Omega\times\Omega
)_{x}}\widehat{
\bP}_{\eta,\xi}(\tau<\infty)
\\[-2pt]
&&\qquad \leq\sum_{x\in\bZ^d}\sup_{(\eta,\xi)\in(\Omega\times
\Omega)_{x}}\int
_0^\infty\sum_{y\in\bZ^d}
w_\alpha(y) \widehat{\bE}_{\eta,\xi}\rho\bigl(\eta_t^1
\bigl(y+X_t^1\bigr), \eta_t^2
\bigl(y+X_t^1\bigr)\bigr) \,dt,\vadjust{\goodbreak}
\end{eqnarray*}
and Lemma~\ref{lemmatriplenorm-coupling-estimate} completes the proof,
where $C_b$ is the combination of the various factors in front of
$\interleave f \interleave$.
\end{pf*}
\begin{pf*}{Proof of Theorem~\ref{thmmu-continuous}}
Let $\alpha, \alpha'$ be two different transition rates. The goal is
to show that
\[
\bigl\llvert \mu^{\mathrm{EP}}_{\alpha}(f)-\mu^{\mathrm{EP}}_{\alpha'}(f)
\bigr\rrvert \leq C \interleave f \interleave
\]
for all $f\dvtx \Omega\to\bR$ with $\interleave f \interleave<\infty$.

The idea is now to use a coupling $\widehat{\bP}$ similar to the one in
Proposition~\ref{propcoupling-construction}.
The coupling contains as objects two copies of the environment, $\eta
^1$ and $\eta^2$, and three random walks, $X^1,X^{12}$ and $X^2$. The
random walk $X^1$ moves on the environment $\eta^1$ with rates $\alpha
$, and correspondingly the random walk $X^2$ moves on $\eta^2$ with
rates~$\alpha'$. The mixed walker $X^{12}$ moves on the environment
$\eta^2$ as well, but according to the rates $\alpha$. The walkers
$X^1,X^2$ will perform the same jumps as $X^{12}$ with maximal
probability. This can be achieved with the same construction as in
Proposition~\ref{propcoupling-construction}, but with Poisson clocks
$N^z$ which have rates $\lambda_z= \sup_{\eta\in\Omega}\alpha
(\eta,z)\vee\alpha'(\eta,z)$.

We only consider the case where all three walkers start at the origin.
We denote by $S_t^{\mathrm{EP},1}$, $S_t^{\mathrm{EP},2}$ the semigroups of the
environment process which correspond to the rates $\alpha$ and $\alpha'$.
Let $\tau:=\inf\{t\geq0\dvtx  X^1_t \neq X^{12}_t \mbox{ or } X^{12}_t
\neq X^{2}_t\}$.
\begin{eqnarray*}
&&
S_t^{\mathrm{EP},1}f(\eta) - S_t^{\mathrm{EP},2}f(\xi)
\\
&&\qquad = \widehat{\bE}_{\eta,\xi} \bigl( f\bigl(\theta_{-X_t^1}
\eta^1_t\bigr) - f\bigl(\theta_{-X_t^2}
\eta^2_t\bigr) \bigr)
\\
&&\qquad = \widehat{\bE}_{\eta,\xi} \ind_{\tau>t} \bigl( f\bigl(\theta
_{-X_t^1}\eta^1_t\bigr) - f\bigl(
\theta_{-X_t^1}\eta^2_t\bigr) \bigr)
\\
&&\qquad\quad{} + \widehat{\bE}_{\eta,\xi} \ind_{\tau\leq t} \bigl( f\bigl(\theta
_{-X_t^1}\eta^1_t\bigr) - f\bigl(
\theta_{-X_t^2}\eta^2_t\bigr) \bigr)
\\
&&\qquad = \widehat{\bE}_{\eta,\xi} \ind_{\tau>t} \bigl( f\bigl(\theta
_{-X_t^1}\eta^1_t\bigr) - f\bigl(
\theta_{-X_t^1}\eta^2_t\bigr) \bigr)
\\
&&\qquad\quad{} + \widehat{\bE}_{\eta,\xi} \ind_{\tau\leq t} \bigl(
S_{t-\tau
}^{\mathrm{EP},1}f\bigl(\theta_{-X_\tau^1}
\eta^1_\tau\bigr) - S_{t-\tau
}^{\mathrm{EP},2}f\bigl(
\theta_{-X_\tau^2}\eta^2_\tau\bigr) \bigr).
\end{eqnarray*}

Therefore,
%
\begin{eqnarray}\label{eqPsiT}
\Psi(T):\!&=&\sup_{0\leq T'\leq T}\sup_{\eta,\xi\in\Omega} \int
_0^{T'} S_t^{\mathrm{EP},1}f(\eta) -
S_t^{\mathrm{EP},2}f(\xi) \,dt
\nonumber
\\
&\leq&\sup_{0\leq T'\leq T}\sup_{\eta,\xi\in\Omega} \int
_0^{T'} \widehat{\bE}_{\eta,\xi}
\ind_{\tau>t} \bigl( f\bigl(\theta_{-X_t^1}\eta^1_t
\bigr) - f\bigl(\theta_{-X_t^1}\eta^2_t\bigr)
\bigr)
\nonumber
\\
&&\hspace*{75.2pt}{} + \widehat{\bE}_{\eta,\xi} \ind_{\tau\leq t}\sup_{\eta,\xi\in
\Omega}
\bigl( S_{t-\tau}^{\mathrm{EP},1}f(\eta) - S_{t-\tau}^{\mathrm{EP},2}f(
\xi) \bigr)\,dt
\\
& \leq&\sup_{0\leq T'\leq T}\sup_{\eta,\xi\in\Omega} \biggl(\widehat{
\bE} _{\eta,\xi} \int_0^\tau f\bigl(
\theta_{-X_t^1}\eta^1_t\bigr) - f\bigl(\theta
_{-X_t^1}\eta^2_t\bigr) \,dt +
\ind_{\tau\leq T'}\Psi\bigl(T'-\tau\bigr) \biggr)
\nonumber
\\
& \leq&\sup_{\eta,\xi\in\Omega} \widehat{\bE}_{\eta,\xi} \biggl(\int
_0^\infty f\bigl(\theta_{-X_t^1}
\eta^1_t\bigr) - f\bigl(\theta_{-X_t^1}
\eta^2_t\bigr) \, dt + \ind_{\tau\leq T}\Psi(T-\tau)
\biggr).\nonumber
\end{eqnarray}
We will now exploit this recursive bound on $\Psi$.
%
\begin{lemma}
Let $\tau_1:= \inf\{t\geq0\dvtx  X^1_t \neq X^{12}_t\}$ and $\tau_2:=
\inf\{t\geq0\dvtx  X^{12}_t \neq X^{2}_t\}$.
Set
\begin{eqnarray*}
\beta&:=&\sum_{z\in\bZ^d} \sup_{\eta\in\bZ^d}
\bigl\llvert \alpha (\eta,z)-\alpha'(\eta,z) \bigr\rrvert,
\\
p(\alpha)&:=&\inf_{\eta,\xi\in\Omega} \widehat{\bP}_{\eta,\xi}(\tau
_1=\infty),
\\
C(\alpha) &:=& \int_0^\infty \bigl( \bigl\llVert
\gamma^+(\alpha )-\gamma^-(\alpha) \bigr\rrVert _\infty t+1
\bigr)^d\sup_{\eta,\xi\in
\Omega} \widehat{\bE}
{}^E_{\eta,\xi}\rho\bigl(\eta^1_t(0),
\eta^2_t(0)\bigr)\,dt,
\end{eqnarray*}
where $\gamma^+(\alpha),\gamma^-(\alpha)$ are as in Proposition
\ref{propcoupling-construction} for the rates $\alpha$.

Let $Y\in\{0,1\}$ be Bernoulli with parameter $p(\alpha)$ and $Y'$
exponentially distributed with parameter $\beta$. Let $Y_1,Y_2,\ldots$ be
i.i.d. copies of $Y\cdot Y'$ and
$N(T):=\inf\{N\geq0\dvtx  \sum_{n=1}^N Y_n>T\}$. Then
\[
\Psi(T) \leq C(\alpha)\interleave f \interleave\bE N(T).
\]
\end{lemma}
\begin{pf}
By construction of the coupling, $\tau_2$ stochastically dominates
$Y'$. As we have
$\tau=\tau_1\wedge\tau_2$ it follows that $\tau\succeq Y_1$.
Using this fact together with the monotonicity of $\Psi$ in (\ref{eqPsiT}),
\begin{eqnarray*}
\Psi(T) &\leq& \sup_{\eta,\xi\in\Omega} \biggl(\widehat{\bE}_{\eta,\xi}
\int_0^\infty f\bigl(\theta_{-X_t^1}
\eta^1_t\bigr) - f\bigl(\theta_{-X_t^1}\eta
^2_t\bigr) \,dt + \ind_{\tau\leq T}\Psi(T-\tau)
\biggr)
\\
&\leq& \sup_{\eta,\xi\in\Omega} \widehat{\bE}_{\eta,\xi} \int
_0^\infty f\bigl(\theta_{-X_t^1}
\eta^1_t\bigr) - f\bigl(\theta_{-X_t^1}
\eta^2_t\bigr) \,dt + \bE \ind_{Y_1\leq T}
\Psi(T-Y_1).
\end{eqnarray*}
As $p(\alpha)>0$ by Lemma~\ref{lemmadecoupling-estimate} we can
iterate this estimate until it terminates after $N(T)$ steps. Therefore
we obtain
\[
\Psi(T) \leq\bE N(T) \sup_{\eta,\xi\in\Omega} \widehat{\bE}_{\eta,\xi}
\int_0^\infty f\bigl(\theta_{-X_t^1}
\eta^1_t\bigr) - f\bigl(\theta_{-X_t^1}\eta
^2_t\bigr) \,dt.
\]
The integral is estimated by telescoping over single site changes and
Lem\-ma~\ref{lemmasupnorm-coupling-estimate} in the usual way, yielding
\[
\Psi(T) \leq C(\alpha)\interleave f \interleave\bE N(T).
\]
\upqed\end{pf}

To finally come back to the original question of continuity,
\begin{eqnarray*}
\bigl\llvert \mu^{\mathrm{EP}}_{\alpha}(f)-\mu^{\mathrm{EP}}_{\alpha'}(f)
\bigr\rrvert &=& \frac1T \biggl\llvert \iiint_0^T
S_t^{\mathrm{EP},1}f(\eta) - S_t^{\mathrm{EP},2}f(\xi)\,dt
\,\mu^{\mathrm{EP}}_{\alpha}(d\eta) \mu^{\mathrm{EP}}_{\alpha'}(d
\xi) \biggr\rrvert
\\
&\leq&\frac1T \Psi(T) \leq\frac{1}{T} \bE N(T) C(\alpha) \interleave f
\interleave
\konv{T\to\infty} \frac{1}{\bE YY'} C(\alpha) \interleave f \interleave
\\
&=&\frac{C(\alpha)}{p(\alpha)}\sum_{z\in\bZ^d} \sup
_{\eta\in\bZ
^d}\bigl\llvert \alpha(\eta,z)-\alpha'(\eta,z)
\bigr\rrvert \interleave f \interleave.
\end{eqnarray*}
By sending $\alpha'$ to $\alpha$, the right-hand side tends to 0 so
that the ergodic measure of the environment process is indeed
continuous in the rates $\alpha$. It is also interesting to note that
both $p(\alpha)$ and $C(\alpha)$ are rather explicit given the original
coupling of the environment. Notably when $\alpha(\eta,z)=\alpha(z)$,
that is, the rates do not depend on the environment, $p(\alpha)=1$ and
$C(\alpha) = \int_0^\infty\sup_{\eta,\xi\in
\Omega}\widehat{\bE}{}^E_{\eta,\xi}\rho(\eta_t^1(0),\eta_t^2(0))\,dt$.
\end{pf*}
\begin{pf*}{Proof of Theorem~\ref{thmEP-phi-estimate}}
The proof of this theorem is mostly identical to the proof of Theorem
\ref{thmEP-estimate}. Hence instead of copying the proof, we just
state where details differ.

A first fact is that the conditions for (a) and (b) imply Assumptions
\ref{assumption1a} and~\ref{assumption1b}. In the adaptation of the proof for part (a), in most
lines it suffices to add a $\phi (\frac{t}{K} )$ to the
integrals. However, in line
(\ref{eqphi-1}), we use
%
\begin{equation}\label{eqphiK-factor}
\phi \biggl(\frac{t}{K} \biggr)\leq\phi \biggl(\frac{t-\tau
}{K} \biggr)
\phi \biggl(\frac{\tau}{K} \biggr)
\end{equation}
to obtain the estimate
\[
\widehat{\bE}_{\eta,\xi} \phi \biggl(\frac{\tau}{K} \biggr) \int
_0^{(T-\tau
)\vee0} \phi \biggl(\frac{t}{K} \biggr)
\bigl\llvert S^{\mathrm{EP}}_{t}f\bigl(\theta _{-X^1_\tau}
\eta^1_\tau\bigr)-S^{\mathrm{EP}}_{t}f\bigl(
\theta_{-X^2_\tau}\eta ^2_\tau\bigr) \bigr\rrvert \,dt
\]
instead. Thereby in lines (\ref{eqphi-2}), (\ref{eqphi-3}) and (\ref
{eqphi-4}) we have to change $\widehat{\bP}_{\eta,\xi}(\tau<\infty)$ to
$\widehat{\bE}
_{\eta,\xi}\phi (\frac{\tau}{K} )\ind_{\tau<\infty
}$. This change then
leads to the replacement of
\[
\inf_{\eta,\xi\in\Omega}\widehat{\bP}_{\eta,\xi}(\tau=\infty)
\]
by the term
\[
1-\sup_{\eta,\xi\in\Omega}\widehat{\bE}_{\eta,\xi}\phi \biggl(
\frac
{\tau}{K} \biggr)\ind _{\tau<\infty}
\]
in the lines (\ref{eqphi-5}) and (\ref{eqphi-6}) [where naturally
$C_a$ becomes $C_a(K)$].
So all we have to prove that for sufficiently big $K$,
\[
\sup_{\eta,\xi\in\Omega}\widehat{\bE}_{\eta,\xi}\phi \biggl(
\frac{\tau
}{K} \biggr)\ind_{\tau
<\infty} < 1.
\]
In a first step, we show that
\[
\sup_{\eta,\xi\in\Omega}\widehat{\bE}_{\eta,\xi}\phi(\tau)
\ind_{\tau
<\infty} < \infty.
\]

As we already saw in the proof of Lemma \ref
{lemmadecoupling-estimate}, we can view the event of decoupling as the
first jump of a Poisson process with time-dependent and random rates
[equation (\ref{eqdynamic-rate-poisson})]. Hence we have
\begin{eqnarray*}
&&
\widehat{\bE}_{\eta,\xi}\phi(\tau)\ind_{\tau<\infty} \\
&&\qquad= \int
_0^\infty\phi(t)\,d\widehat{\bP}_{\eta,\xi}(
\tau>t)
\\
&&\qquad= \int_0^\infty\phi(t) \widehat{
\bE}_{\eta,\xi} \sum_{z\in\bZ^d}\bigl\llvert \alpha
\bigl(\theta_{-X^1_t}\eta^1_t,z\bigr)- \alpha\bigl(
\theta _{-X^1_t}\eta^2_t,z\bigr) \bigr\rrvert
\\
&&\qquad\quad\hspace*{17.2pt}{}\times\exp \biggl({-\int_0^t \sum
_{z\in\bZ^d}\bigl\llvert \alpha \bigl(\theta_{-X^1_s}
\eta^1_s,z\bigr)- \alpha\bigl(\theta_{-X^1_s}
\eta^2_s,z\bigr) \bigr\rrvert \,ds} \biggr) \,dt
\\
&&\qquad\leq\int_0^\infty\phi(t) \widehat{
\bE}_{\eta,\xi} \sum_{z\in\bZ
^d}\bigl\llvert \alpha
\bigl(\theta_{-X^1_t}\eta^1_t,z\bigr)- \alpha\bigl(
\theta _{-X^1_t}\eta^2_t,z\bigr) \bigr\rrvert \,dt.
\end{eqnarray*}
By telescoping over single site discrepancies and using Lemma \ref
{lemmasupnorm-coupling-estimate}, this is less than
\[
\int_0^\infty\phi(t) \bigl(\bigl\llVert \gamma^+-
\gamma^- \bigr\rrVert _\infty+1\bigr)^d t^d \sup
_{\eta,\xi\in\Omega} \widehat{\bE}{}^E_{\eta,\xi}\rho\bigl(\eta
_t^1(0),\eta _t^2(0)\bigr) \,dt<
\infty
\]
by assumption.
Since $\phi(t/K)$ decreases to 1 as $K\to\infty$, monotone
convergence implies
\[
\lim_{K\to\infty} \widehat{\bE}_{\eta,\xi}\phi \biggl(
\frac{\tau
}{K} \biggr)\ind_{\tau<\infty
} = \widehat{\bE}_{\eta,\xi}
\ind_{\tau<\infty}< 1
\]
by Lemma~\ref{lemmadecoupling-estimate}. Consequently, there exists a
$K_0\geq0$ such that for all $K>K_0$,
\[
\widehat{\bE}_{\eta,\xi}\phi \biggl(\frac{\tau}{K} \biggr)
\ind_{\tau
<\infty} < 1.
\]
This completes the adaptation of part (a).

The adaptation of the proof of part (b) follows the same scheme, where
we add the term $\phi (\frac{t}{K} )$ to all integrals. Note
that this gives a version of~Lem\-ma~\ref{lemmatriplenorm-coupling-estimate} as well. Then, in line
(\ref{eqphi-7}) we use (\ref{eqphiK-factor}) again and then have to
replace $\widehat{\bP}_{\eta,\xi}(\tau<\infty)$ by
$\widehat{\bE}_{\eta,\xi}\phi (\frac{\tau}{K} )\ind_{\tau
<\infty}$
in lines (\ref{eqphi-8}) and (\ref{eqphi-9}). To estimate
(\ref{eqphi-9}), we use\looseness=-1
\begin{eqnarray*}
&&
\widehat{\bE}_{\eta,\xi}\phi \biggl(\frac{\tau}{K} \biggr)
\ind_{\tau
<\infty} \\
&&\qquad\leq\int_0^\infty\phi \biggl(
\frac{t}{K} \biggr) \widehat{\bE}_{\eta,\xi} \sum
_{z\in\bZ
^d}\bigl\llvert \alpha\bigl(\theta_{-X^1_t}
\eta^1_t,z\bigr)- \alpha\bigl(\theta _{-X^1_t}
\eta^2_t,z\bigr) \bigr\rrvert \,dt
\\
&&\qquad\leq\int_0^\infty\phi \biggl(\frac{t}{K}
\biggr) \sum_{y\in\bZ
^d} w_{\alpha} \widehat{\bE}
_{\eta,\xi} \rho\bigl(\eta^1_t\bigl(y+X_t^1
\bigr),\eta^2_t\bigl(y+X_t^1\bigr)
\bigr) \,dt
\end{eqnarray*}\looseness=0
with $w_\alpha$ as in the original proof. Therefore
\begin{eqnarray*}
&&
\sum_{x\in\bZ^d}\sup_{(\eta,\xi)\in(\Omega\times\Omega
)_{x}}\widehat{
\bE}_{\eta,\xi}\phi \biggl(\frac{\tau}{K} \biggr)\ind_{\tau
<\infty}
\\
&&\qquad \leq\sum_{x\in\bZ^d}\sup_{(\eta,\xi)\in(\Omega\times
\Omega)_{x}}\int
_0^\infty\phi \biggl(\frac{t}{K} \biggr)\sum
_{y\in\bZ^d} w_\alpha (y) \widehat{
\bE}_{\eta,\xi}\rho\bigl(\eta_t^1\bigl(y+X_t^1
\bigr),\\
&&\qquad\quad\hspace*{214pt} \eta_t^2\bigl(y+X_t^1
\bigr)\bigr) \, dt,
\end{eqnarray*}
which is finite by Lemma
\ref{lemmatriplenorm-coupling-estimate}.
\end{pf*}
\begin{pf*}{Proof of Proposition~\ref{propLLN}}
The LLN for $\bP_{\mu^\mathrm{EP},0}$ follows directly by ergodicity. To prove
the LLN for $\bP_{\eta,x}$ we use a slight modification of the
coupling $\widehat{\bP}$ in Proposition~\ref{propcoupling-construction}.

Suppose w.l.o.g. that $x=0$ (otherwise look at
$\bP_{\theta_{-x}\eta,0}$). In the construction of the modified
coupling $\widetilde\bP$ we look at $X^1_T-X^2_T$. Jump events of
$X^1_T-X^2_T$ we call decoupling events, which are events when one
walker jumps but the other does not. Up to the first decoupling the
coupling $\widetilde\bP_{\eta,0;\mu^{\mathrm{EP}},0}$ is identical to
$\widehat{\bP}_{\eta,0;\mu^{\mathrm{EP}},0}$. By Lemma \ref
{lemmadecoupling-estimate}, there is at least probability $p>0$
uniformly in $\eta$ to never\vspace*{1pt} decouple. At the instant $\tau$ of a
decoupling event we restart the coupling $\widehat{\bP}{}^E$ of the
environment. This is done in the configuration
$\theta_{-X^1_\tau}\eta^{1}_{\tau},\theta_{-X^2_\tau}\eta^{2}_{\tau}$.
That is, instead of coupling $\eta^1_t(x)$ with $\eta^2_t(x)$, we match
$\eta^1_t(x+X^1_\tau)$ with $\eta^2_t(x+X^2_\tau)$. This allows us to
apply Lemma~\ref{lemmadecoupling-estimate} a second time, since for
the purpose of decoupling events, both walkers start at the origin at
time $\tau$. Iterating, we then have at most a geometric number $N$ of
decoupling events at $\tau_1,\ldots,\tau_N$. Hence
\[
\bigl\llvert X^1_T-X^2_T \bigr
\rrvert \leq\sum_{n=1}^N \bigl(\bigl
\llvert X^1_{\tau_n}-X^1_{\tau_n-} \bigr
\rrvert +\bigl\llvert X^2_{\tau
_n}-X^2_{\tau_n-}
\bigr\rrvert \bigr),
\]
which when divided by $T$ converges to 0 in $L^1$ and almost surely
w.r.t. $\widetilde\bP_{\eta,0;\mu^{\mathrm{EP}},0}$.

Therefore $X_T/T$ converges in $L^1(\bP_{\eta,x})$, for given $x$
uniformly in $\eta\in\Omega$, and $\bP_{\eta,x}$ almost surely to
the same limit, $\lim_{T\to\infty}\frac{1}{T}\bE_{\mu^{\mathrm{EP}},0}(X_T)$.
\end{pf*}



%

\printaddresses

\end{document}